\documentclass[10pt, onecolumn]{IEEEtran}

\pdfoutput=1

\usepackage{algorithm}
\usepackage{algpseudocode}

\usepackage{amsmath}
\usepackage{algorithm}
\usepackage{array}
\usepackage{mdwmath}
\usepackage{eqparbox}                                                                                                                           
\usepackage{float}
\usepackage[font=footnotesize]{subfig}
\usepackage{fixltx2e}
\usepackage{stfloats}
\usepackage{amsfonts}
\usepackage{cite}
\usepackage{amsmath}
\usepackage{bbm}
\usepackage{amssymb}
\usepackage{mathrsfs}
\usepackage{amsbsy}	
\usepackage{color}
\usepackage{multicol}

\usepackage{graphicx,epstopdf}

\usepackage{verbatim}

\floatstyle{ruled}
\newfloat{algorithm}{tbp}{loa}
\floatname{algorithm}{Algorithm}

\long\def\hide#1{}

\usepackage{comment}

\newtheorem{lemma}{\bf Lemma}[section]

\newtheorem{proposition}{\bf Proposition}[section]

\newcommand{\droped}[1]{{\color{blue} \sout{}}}

\def\ba{\begin{array}}
\def\ea{\end{array}}
\newcommand{\beq}{\begin{equation}}
\newcommand{\eeq}{\end{equation}}
\newcommand{\bq}{\begin{eqnarray}}
\newcommand{\eq}{\end{eqnarray}}
\newcommand{\bqn}{\begin{eqnarray*}}
\newcommand{\eqn}{\end{eqnarray*}}
\newcommand{\bee}{\begin{enumerate}}
\newcommand{\eee}{\end{enumerate}}
\newcommand{\bi}{\begin{itemize}}
\newcommand{\ei}{\end{itemize}}

\newcommand{\mathN}{\mathcal{N}}
\newcommand{\mathE}{\mathcal{E}}

\newcommand{\mathT}{\mathcal{T}}
\newcommand{\mathI}{\mathcal{I}}
\newcommand{\mathK}{\mathcal{K}}

\newboolean{showcomments}
\setboolean{showcomments}{true}
\newcommand{\qiuyu}[1]{  \ifthenelse{\boolean{showcomments}}
{ \textcolor{blue}{(Qiuyu says:  #1)}} {}  }
\newcommand{\anwar}[1]{\ifthenelse{\boolean{showcomments}}
{ \textcolor{brown}{(Anwar says: #1)} } {} }
\newcommand{\slow}[1]{\ifthenelse{\boolean{showcomments}}
{ \textcolor{red}{(Steven says:  #1)}}{}}
\newcommand{\pmark}[1]{\ifthenelse{\boolean{showcomments}}
{ \textcolor{cyan}{(Label:  #1)}}{}}

\begin{document}

\title{ Distributed Optimal Power Flow Algorithm for Balanced Radial Distribution Networks}

\author{Qiuyu Peng and Steven H. Low
\thanks{A preliminary version has appeared in \cite{peng2014cdc}.}
\thanks{Qiuyu Peng is with the Electrical Engineering Department 
and Steven H. Low is with the Computing and Mathematical Sciences and
the Electrical Engineering Departments, California Institute of Technology, Pasadena, CA 91125, USA.
{\small \tt \{qpeng, slow\}@caltech.edu}}%
}

\maketitle
\begin{abstract}

The optimal power flow (OPF) problem is fundamental in power system
operations and planning. Large-scale renewable penetration in distribution networks calls for real-time 
feedback control, and hence the need for fast and distributed solutions for OPF. 
This is difficult because OPF is nonconvex and Kirchhoff's laws are global.   In
this paper we propose a solution for \emph{balanced} radial distribution networks.  It exploits recent
results that suggest solving for a globally optimal solution of OPF over a
radial network through the second-order cone program (SOCP) relaxation.
Our distributed algorithm is based on alternating direction method of 
multiplier (ADMM), but unlike standard ADMM algorithms that often require
iteratively solving optimization subproblems in each ADMM iteration, our
decomposition allows us to derive closed form solutions for these subproblems,
greatly speeding up each ADMM iteration.  We present simulations on a real-world
2,065-bus distribution network to illustrate the scalability and optimality of
the proposed algorithm.


\end{abstract}

\begin{IEEEkeywords}
Power Distribution, Nonlinear systems, Optimal Power Flow, Distributed Algorithm
\end{IEEEkeywords}

\section{Introduction}

\IEEEPARstart{T}HE optimal power flow (OPF) problem seeks to optimize certain objective such as power loss and 
generation cost subject to power flow equations and operational constraints. 
It is a fundamental problem because it underlies many power system operations and planning problems
such as economic dispatch, unit commitment, state estimation, stability and reliability assessment, 
volt/var control, demand response, etc.  The continued growth of highly volatile renewable 
sources on distribution systems calls for real-time feedback control.   Solving the OPF problems in such an
environment has at least two challenges.

First the OPF problem is hard to solve because of its nonconvex feasible set.  
Recently a new approach through convex relaxation has been developed.
Specifically semidefinite program (SDP) relaxation \cite{Bai2008} and second order cone program (SOCP) relaxation \cite{Jabr2006} have been proposed in the bus injection model, and SOCP relaxation has been proposed
in the branch flow model \cite{Farivar2011-VAR-SGC, Farivar-2013-BFM-TPS}.  
See the tutorial \cite{Low2014a, Low2014b} for further pointers to the literature.
When an optimal solution of the original OPF problem can be recovered from any optimal solution 
of a convex relaxation, we say the relaxation is \emph{exact}.  
For radial distribution networks (whose graphs are trees), several sufficient conditions have been proved that 
guarantee SOCP and SDP relaxations are exact.  This is important because almost all distribution
systems are radial.  Moreover some of these conditions have been shown to hold for many practical networks.  
In those cases we can rely on off-the-shelf convex optimization solvers to obtain a globally 
optimal solution for the nonconvex OPF problem.

Second most algorithms proposed in the literature are centralized and meant for applications
in today's energy management systems that, \emph{e.g.}, centrally schedule a relatively small number of generators.   
In future networks that simultaneously optimize (possibly real-time) the
operation of a large number of intelligent endpoints, a centralized approach will not scale because of 
 its computation and communication overhead.
 In this paper we address this challenge.  Specifically 
we propose a distributed algorithm for solving the SOCP relaxation of OPF for  balanced radial distribution
networks.

Various distributed algorithms have been developed to solve the OPF problem. Through optimization decomposition, the original OPF problem is decomposed into several local subproblems that can be solved simultaneously. Some distributed algorithms do not deal with the non convexity issue of the OPF problem, including \cite{kim1997coarse,baldick1999fast}, which leverage method of multipliers and \cite{sun2013fully}, which is based on alternating direction method of multiplier (ADMM). However, the convergence of these algorithms is not guaranteed due to non-convexity of the problem. In contrast, algorithms for the convexified OPF problem are proposed to guarantee convergence, \emph{e.g.} dual decomposition method \cite{lam2012distributed,lam2012optimal}, auxiliary variable method  \cite{devane2013stability,li2012demand} and ADMM \cite{dall2012distributed,kraning2013dynamic}.

One of the key performance metrics of a distributed algorithm is the time of convergence (ToC), which is the product of the number of iterations and the computation time to solve the subproblems in each iteration. To our knowledge, all the distributed OPF algorithms in the literature rely on generic iterative optimization solvers, which are computationally intensive, to solve the optimization subproblems. In this paper, we will improve ToC by reducing the computation time for each subproblem.

Specifically we develop a scalable distributed algorithm 
through decomposing the \emph{convexified} OPF problem into smaller subproblems based on alternating direction method of multiplier (ADMM). ADMM blends the decomposability of dual decomposition and superior convergence properties of the method of multipliers \cite{boyd2011distributed}. It has broad applications in different areas and particularly useful when the subproblems can be solved efficiently \cite{ghadimi2013optimal}, for example when they admit closed form expressions, \emph{e.g.} matrix factorization \cite{sun2014alternating}, image recovery \cite{afonso2010fast}.

The proposed algorithm has two advantages: 1) There is closed form solution for each optimization
subproblem, thus eliminating the need for an iterative procedure to solve a SDP/SOCP problem for each ADMM iteration.
2) Communication is only required between adjacent buses. 

We demonstrate the scalability of the proposed algorithms using a real-life network. 
In particular, we show that the algorithm converges within 0.6s for a 2,065-bus system. 
To show the superiority of deriving close form expression of each subproblems, 
finally we compare the computation time for solving a subproblem by our algorithm 
and an off-the-shelf optimization solver (CVX, \cite{grant2008cvx}). 
Our solver requires on average $6.8\times 10^{-4}$s while CVX requires on average $0.5$s.
On the other hand, we also show that the convergence rate is mainly determined by the diameter \footnote{The diameter of a graph is defined as the number of hops between two furthest nodes.} of the network through simulating the algorithm on different networks.

The rest of the paper is structured as follows. The OPF problem is defined in section \ref{sec:opf}. 
In section \ref{sec:alg}, we develop our distributed algorithm.   In section \ref{sec:case}, we
 test its scalability using data from a real-world distribution network. We conclude this paper in section \ref{sec:conclusion}.

%
%
%

\section{Problem formulation}\label{sec:opf}

In this section, we define the optimal power flow (OPF) problem 
on a balanced radial distribution network and review how to solve it through SOCP relaxation.

We denote the set of complex numbers with $\mathbb{C}$, the set of $n$-dimensional complex numbers with $\mathbb{C}^n$. The hermitian transpose of a vector is denoted by $()^H$. To differentiate vector and scaler operations, the conjugate of a complex scaler is denoted by $()^*$. The inner product of two vectors $x,y\in\mathbb{C}^n$ is denoted by $\langle x,y\rangle:=\mathbf{Re}(tr(x^Hy))$. The Euclidean norm of a vector $x\in\mathbb{C}^{n}$ is defined as $\|x\|_2:=\sqrt{\langle x,x\rangle}$.

\subsection{Branch flow model}

\begin{figure}
\centering
\includegraphics[width=0.5\linewidth]{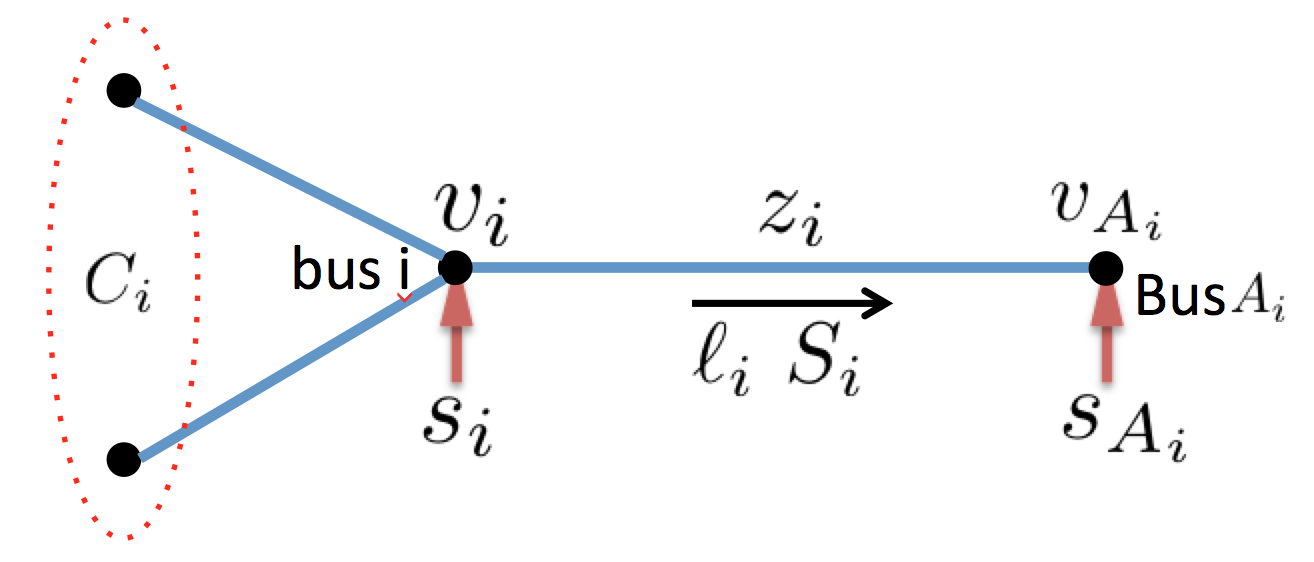}
\caption{Notations.}
\label{fig:notation}
\end{figure}

We model a distribution network by a \emph{directed} tree graph $\mathT := (\mathN, \mathE)$ where
 $\mathN:=\{0,\ldots,n\}$ represents the set of buses and $\mathE$ represents the set
 of distribution lines connecting the buses in $\mathN$.  
Index the root of the tree by $0$ and let $\mathN_+:=\mathN\setminus\{0\}$ denote the other buses. For each node $i$, it has a unique ancestor $A_i$ and a set of children nodes, denoted by $C_i$.
We adopt the graph orientation where every line points towards the root. Each directed line connects a node $i$ and its unique ancestor $A_i$. We hence label the lines by $\mathE:=\{1,\ldots,n\}$ where each 
 $i\in\mathE$ denotes a line from $i$ to $A_i$. 

For each bus $i\in \mathN$, let $V_i=|V_i|e^{\mathbf{i}\theta_i}$ be its complex voltage and $v_i:=|V_i|^2$ be its magnitude squared. Let $s_i:=p_i+\mathbf{i}q_i$ be its net complex power injection which is generation minus load. For each line $i\in \mathE$, let $z_{i}=r_{i}+\mathbf{i}x_{i}$ be its complex impedance. Let $I_{i}$ be the complex branch current from bus $i$ to $A_i$ and $\ell_{i}:=|I_{i}|^2$ be its magnitude squared. Let $S_{i}:=P_{i}+\mathbf{i}Q_{i}$ be the branch power flow from bus $i$ to $A_i$. The notations are illustrated in Fig. \ref{fig:notation}. A variable without a subscript denotes a column vector with appropriate components, as summarized below.
\begin{center}
\begin{tabular}{|c|c|c|}
\hline
 $v:=(v_{i},i\in \mathN)$ &  $s:=(s_{i},i\in \mathN) $\\
 \hline
  $\ell:=(\ell_i,i\in \mathE)$ & $S:=(S_{i},i\in \mathE)$  \\
  \hline
\end{tabular}
\vspace{0.05in}
\end{center}

Branch flow model is first proposed in \cite{Baran1989a,Baran1989b} for radial network. It has better numerical stability than bus injection model and has been advocated for the design and operation for radial distribution network, \cite{Farivar-2013-BFM-TPS, gan2014exact,li2012demand}. It ignores the phase angles of voltages and currents and uses only the set of variables $(v,s,\ell,S)$.  
Given a radial network $\mathT$, the branch flow model is defined by:
\begin{subequations}\label{eq:bfm}
\begin{align}
&v_{A_i}-v_i+(z_iS_i^*+S_iz_i^*)-\ell_{i}|z_i|^2=0 & i\in\mathE  \label{eq:bfm1}\\  
& \sum_{j\in C_i}(S_{j}-\ell_{j}z_{j})+s_i-S_{i}=0  & i\in\mathN  \label{eq:bfm2}\\  
& |S_i|^2= v_i\ell_i  & i\in\mathE   \label{eq:bfm3}
\end{align}
\end{subequations}
where  $S_0=0$ ( the root of the tree does not have parent) for ease of presentation. Given a vector $(v,s,\ell,S)$ that satisfies \eqref{eq:bfm}, the phase angles of the voltages and currents can be uniquely determined if the network is a tree.
Hence the branch flow model  \eqref{eq:bfm} is equivalent to a full  AC power flow model. 
See \cite[Section III-A]{Farivar-2013-BFM-TPS} for details. 

\subsection{OPF and SOCP relaxation}

The OPF problem seeks to optimize certain objective, e.g. total line loss or total generation cost, subject to power flow equations \eqref{eq:bfm} and various operational constraints. 
We consider an objective function of the following form: 
\bq\label{eq:objective}
F(p)=\sum_{i\in\mathN}f_i(p_i):=\sum_{i\in\mathN}\left(\frac{\alpha_i}{2}p_i^2+\beta_ip_i\right)
\eq
where $\alpha_i,\beta_i\geq0$.
For instance,
\bi
\item to minimize total line loss, we can set $\alpha_i=0$ and $\beta_i=1$ for each bus $i\in\mathN$.
\item to minimize generation cost, we can set $\alpha_i=0$ and $\beta_i=0$ for bus $i$ where there is no generator and for generator bus $i$, the corresponding $\alpha_i,\beta_i$ depends on the characteristic of the generator. 
\ei

We consider two operational constraints. First, the power injection $s_i$ at each bus $i$ is constrained to be in a region $\mathI_i$, \emph{i.e.} 
\bq\label{eq:operation1}
s_i\in\mathI_i \ \ \text{for } i\in\mathN
\eq
The feasible power injection region $\mathI_i$ is determined by the controllable devices attached to bus $i$. Some common controllable loads are:
\bi
\item For controllable load, whose real power can vary within $[\underline p_i,\overline p_i]$ and reactive power can vary within $[\underline q_i,\overline q_i]$, the injection region $\mathI_i$ is 
\begin{subequations}
\bq\label{eq:S2}
\mathI_i=\{p+\mathbf{i}q\mid p\in[ \underline p_i,  \overline p_i], q\in[\underline q_i,\overline q_i] \}\subseteq \mathbb{C}
\eq
\item For solar panel connecting the grid through a inverter with nameplate $\overline s_i$, the injection region $\mathI_i$ is 
\bq\label{eq:S1}
\mathI_i=\{p+\mathbf{i}q\mid p\geq 0, p^2+q^2\leq \overline s_i^2\}\subseteq\mathbb{C}
\eq
\end{subequations}
\ei
Second, the voltage magnitude at each bus $i\in\mathN$ needs to be maintained within a prescribed region, \emph{i.e.} 
\bq\label{eq:operation2}
\underline v_i \leq v_i\leq \overline v_i  \ \ \text{for } i\in\mathN
\eq
Typically the voltage magnitude at the substation bus $0$ is assumed to be fixed at some prescribed value, \emph{i.e.} $\underline v_0=\overline v_0$. At other bus $i\in\mathN_+$, the voltage magnitude is typically allowed to deviate by $5\%$ from its nominal value $1$, \emph{i.e.} $\underline v_i=0.95^2$ and $\overline v_i=1.05^2$.

To summarize, the OPF problem for radial network is
\bq
\text{\bf OPF: }\min && \sum_{i\in\mathN}f_i(p_i)\nonumber\\
\mathrm{over} && v,s,S,\ell \label{eq:opf}\\
\mathrm{s.t.} &&  \eqref{eq:bfm}, \eqref{eq:operation1}\text{ and } \eqref{eq:operation2}\nonumber
\eq

The OPF problem \eqref{eq:opf} is nonconvex due to the quadratic equality constraint \eqref{eq:bfm3}. In \cite{Farivar2011-VAR-SGC, Farivar-2013-BFM-TPS}, \eqref{eq:bfm3} is relaxed to a second order cone constraint:
\bq\label{eq:socp}
|S_i|^2\leq v_i\ell_i  \quad\text{for } \quad i\in\mathE,
\eq
resulting in a second-order cone program (SOCP) relaxation of \eqref{eq:opf} 
\bq
\text{\bf ROPF: }\min && \sum_{i\in\mathN}f_i(p_i)\nonumber\\
\mathrm{over} && v,s,S,\ell \label{eq:ropf}\\
\mathrm{s.t.} &&  \eqref{eq:bfm1}, \ \eqref{eq:bfm2}, \ \eqref{eq:socp}  \text{ and }\eqref{eq:operation1},\eqref{eq:operation2}\nonumber
\eq

Clearly the relaxation ROPF \eqref{eq:ropf} provides a lower bound for the original OPF problem \eqref{eq:opf} since the original feasible set is enlarged. The relaxation is called \emph{exact} if every optimal solution of ROPF attains equality in \eqref{eq:bfm3} and hence is also optimal for the original OPF. 
For network with tree topology, SOCP relaxation is exact under some mild conditions 
\cite{Farivar-2013-BFM-TPS, gan2014exact}.

\section{Distributed Algorithm for OPF}\label{sec:alg}

We assume SOCP relaxation is exact and  develop in this section
a distributed algorithm that solves ROPF. 
We first review a standard alternating direction method of multiplier (ADMM).
We then make use of the structure of ROPF to speed up the standard ADMM algorithm 
by deriving closed form expressions for the optimization subproblems in each ADMM iteration.

\subsection{Preliminary: ADMM}

ADMM blends the decomposability of dual decomposition with the superior convergence properties of the method of multipliers
\cite{boyd2011distributed}. For our application, we consider optimization problems of the form:
\bq
\min && f(x)+g(z) \nonumber \\
\mathrm{over} &&  x\in\mathcal{K}_x, \ \ z\in\mathcal{K}_z \label{eq:admm}\\
\mathrm{s.t.} &&  x=z \nonumber
\eq
where $\mathcal{K}_x,\mathcal{K}_z$ are convex sets. Let $\lambda$ denote the Lagrange
 multiplier for the constraint $x=z$.   Then the augmented Lagrangian is defined as
\bq\label{eq:agumentlag}
L_\rho(x,z,\lambda):=f(x)+g(z)+\langle \lambda, x-z\rangle+\frac{\rho}{2}\|x-z\|_2^2,
\eq
where  $\rho\geq 0$ is a constant. When $\rho=0$, the augmented Lagrangian reduces to 
the standard Lagrangian. At each iteration $k$, ADMM consists of the iterations:
\begin{subequations}\label{eq:update}
\bq
x^{k+1}&\in&\arg\min_{x\in\mathcal{K}_x} L_\rho(x,z^{k},\lambda^k)\label{eq:xupdate}\\
z^{k+1}&\in&\arg\min_{z\in\mathcal{K}_z} L_\rho(x^{k+1},z,\lambda^k)\label{eq:zupdate}\\
\lambda^{k+1}&=&\lambda^{k}+\rho(x^{k+1}-z^{k+1})\label{eq:mupdate}
\eq
\end{subequations}
Compared to dual decomposition, ADMM is guaranteed to converge to an optimal 
solution under less restrictive conditions. Let
\begin{subequations}\label{eq:feasible}
\bq
r^k&:=&\|x^{k}-z^{k}\|_2 \label{eq:pfeasible} \\ 
s^k&:=&\rho\|z^{k}-z^{k-1}\|_2 \label{eq:dfeasible}
\eq
\end{subequations}
which can be viewed as the residuals for primal and dual feasibility.
Assume: 
\bi
\item A1: $f$ and $g$ are closed proper and convex.
\item A2: The unaugmented Lagrangian $L_0$ has a saddle point.
\ei
The correctness of ADMM is guaranteed by the following result; see \cite[Chapter 3]{boyd2011distributed}.
\begin{proposition}[\cite{boyd2011distributed}]\label{prop}
Suppose A1 and A2 hold. Let $p^*$ be the optimal objective value.  
Then 
\bqn
\lim_{k\rightarrow \infty}r^k=0,  \quad  \lim_{k\rightarrow \infty}s^k=0
\eqn
and
\bqn
\lim_{k\rightarrow \infty} f(x^k)+g(z^k)= p^*
\eqn
\end{proposition}

\subsection{Apply ADMM to OPF problem}

We assume the SOCP relaxation is exact and now derive a distributed algorithm for solving ROPF \eqref{eq:ropf} that has the
following advantages: 
\bi
\item Each bus only needs to solves a local subproblem in each iteration of \eqref{eq:update}.  
	Moreover there is a closed form solution for each subproblem, in contrast to most 
	algorithms that employ iterative procedure to solve these subproblems
\cite{lam2012optimal,lam2012distributed,dall2012distributed,sun2013fully,devane2013stability,li2012demand,kim1997coarse,baldick1999fast}.
\item Communication is only required between adjacent buses.
\ei

The ROPF problem defined in \eqref{eq:ropf} can be written explicitly as:
\begin{subequations}
\begin{align}
\min & \ \ \sum_{i\in \mathN} f_i(p_i) &\\
\mathrm{over }& \ \ v,s,S,\ell& \\
\mathrm{s.t.}  & \ \ v_{A_i}-v_i+z_iS_i^*+S_iz_i^*-\ell_i|z_i|^2=0 & i\in\mathE  \label{eq:bfm21}\\
&\ \ \sum_{i\in C_i}(S_j-z_j\ell_j)-S_i +s_i=0 & i\in\mathN  \label{eq:bfm22}\\
&\ \ |S_i|^2\leq v_i\ell_i  &  i\in\mathE  \label{eq:bfm23}\\
&\ \ s_i\in\mathI_i & i\in\mathN   \label{eq:bfm24} \\
&\ \ \underline v_i\leq v_i\leq \overline v_i & i\in\mathN   \label{eq:bfm25}
\end{align}
\end{subequations}

Assume each bus $i$ is an agent that maintains local variables $(v_i,s_i,S_i,\ell_i)$. Then \eqref{eq:bfm23}--\eqref{eq:bfm25} are local constraints to agent (bus) $i$.  \eqref{eq:bfm21} and \eqref{eq:bfm22} describe the coupling constraints among $i$ and its parent $A_i$ and the set of children in $C_i$, \emph{i.e.} \eqref{eq:bfm21} models the voltage of its ancestor $A_i$ as a function of the local variables of $i$, \eqref{eq:bfm22} describes the power flow balance among the set of children $C_i$ and bus $i$ itself. To decouple the constraints  \eqref{eq:bfm21}--\eqref{eq:bfm22}, for each bus $i$, its ancestor $A_i$ sends its voltage $v_{A_i}$ to $i$, denoted by $v_{A_i,i}$ and each child $j\in C_i$ sends the branch power to $i$, denoted by $S_{j,i}$ and current to $i$, denoted by $\ell_{j,i}$. Then ROPF can be written equivalently as follows.

\begin{subequations}\label{eq:eropf}
\begin{align}
\min & \ \ \sum_{i\in \mathN} f_i(p_i^{(z)}) \\
\mathrm{over} & \ \ x:=\{v^{(x)}_i, s^{(x)}_i,S^{(x)}_i,\ell^{(x)}_i,v^{(x)}_{A_i,i},S^{(x)}_{i,A_i},\ell^{(x)}_{i,A_i},i\in\mathN \} \hspace{-0.1in} \nonumber\\
& \ \ z:=\{v^{(z)}_i, s^{(z)}_i,S^{(z)}_i,\ell^{(z)}_{i},i\in\mathN\} \nonumber \\
\mathrm{s.t.}  & \ \ v^{(x)}_{A_i,i}-v^{(x)}_{i}+z_i\left(S^{(x)}_{i}\right)^*+S^{(x)}_{i}z_i^*-\ell^{(x)}_{i}|z_i|^2=0 \hspace{-0.1in}& i\in\mathE  \\
& \sum_{i\in C_i}\left(S^{(x)}_{j,i}-z_j\ell^{(x)}_{j,i}\right)-S^{(x)}_{i} +s^{(x)}_{i}=0 & i\in\mathN \\
& |S_i^{(z)}|^2\leq v_i^{(z)}\ell_i^{(z)} &  i\in\mathE  \\
&  s_i^{(z)}\in \mathI_i  &  i\in\mathN \\
& \underline v_i\leq v_i^{(z)}\leq \overline v_i  &  i\in\mathN \\
& S_{i,A_i}^{(x)}=S_i^{(z)}, \ \ell_{i,A_i}^{(x)}=\ell_{i}^{(z)}, \ v_{A_i,i}^{(x)}=v_{A_i}^{(z)} & i\in\mathN \label{eq:ndup} \\
& S_{i}^{(x)}=S_i^{(z)},  \ell_{i}^{(x)}=\ell_{i}^{(z)}, v_{i}^{(x)}=v_{i}^{(z)}, \ s_i^{(x)}=s_i^{(z)} & i\in\mathN \label{eq:sdup}
\end{align}
\end{subequations}
{\flushleft where \eqref{eq:ndup} and \eqref{eq:sdup} are consensus constraints that force all the copies of each variable to be the same. Since ADMM has two separate groups of variables $x$ and $z$ that is updated alternatively, we put superscripts $(\cdot)^{(x)}$ and $(\cdot)^{(z)}$ on each variable to denote whether the variable is updated in the $x$-update or $z$-update step.}

Next, we apply ADMM to decompose \eqref{eq:eropf} by relaxing the consensus constraints in \eqref{eq:ndup} and \eqref{eq:sdup}. Let $\lambda,\mu,\gamma$ be the Lagrangian multipliers associated with \eqref{eq:ndup} and \eqref{eq:sdup} as specified in Table \ref{tab:muliplier}.

\begin{table}
\caption{Multipliers associated with constraints \eqref{eq:ndup}-\eqref{eq:sdup}}
\begin{center}
\begin{tabular}{|c|c||c|c|}
\hline
$\lambda_{1,i}$: & $S_{i}^{(x)}=S_i^{(z)}$ & $\lambda_{2,i}$: & $\ell_{i}^{(x)}=\ell_{i}^{(z)}$\\
\hline
$\lambda_{3,i}$: & $v_{i}^{(x)}=v_{i}^{(z)}$ & $\lambda_{4,i}$: & $s_i^{(x)}=s_i^{(z)}$\\
\hline
$\mu_{1,i}$: & $S_{i,A_i}^{(x)}=S_i^{(z)}$ & $\mu_{2,i}$: & $\ell_{i,A_i}^{(x)}=\ell_{i}^{(z)}$ \\
\hline
$\gamma_{i}$: & $v_{A_i,i}^{(x)}=v_{A_i}^{(z)}$ & &\\ 
\hline
\end{tabular}
\end{center}
\label{tab:muliplier}
\end{table}

Denote
\bqn
x_i&:=&\left(v_i^{(x)},\ell_i^{(x)},S_i^{(x)},s_i^{(x)}\right)\\
x_{j,i}&:=& \left(\ell_{j,i}^{(x)},S_{j,i}^{(x)}\right)\\
z_i&:=&\left(v_i^{(z)},\ell_i^{(z)},S_i^{(z)},s_i^{(z)}\right)\\
\lambda_i&:=&(\lambda_{k,i},k=1,2,3,4)\\
\mu_i&:=&(\mu_{k,i},k=1,2)
\eqn

The variables maintained by each agent (bus) $i$ are its local variables for itself: $x_i$, $z_i$, the copy of its parent's voltage $v_{A_i,i}^{(x)}$, the copy $x_{j,i}$ from each of its child $j\in C_i$ and the associated Lagrangian multipliers. Let $\mathcal{A}_i$ denote the set of variables, then 
\bqn
\mathcal{A}_i:=\{x_i, v_{A_i,i}^{(x)}, \{x_{j,i},\mu_j,j\in C_i\},z_i,\lambda_i,\gamma_i\}
\eqn

\begin{figure}
\centering
\subfloat[x-update]{
	\includegraphics[scale=0.25]{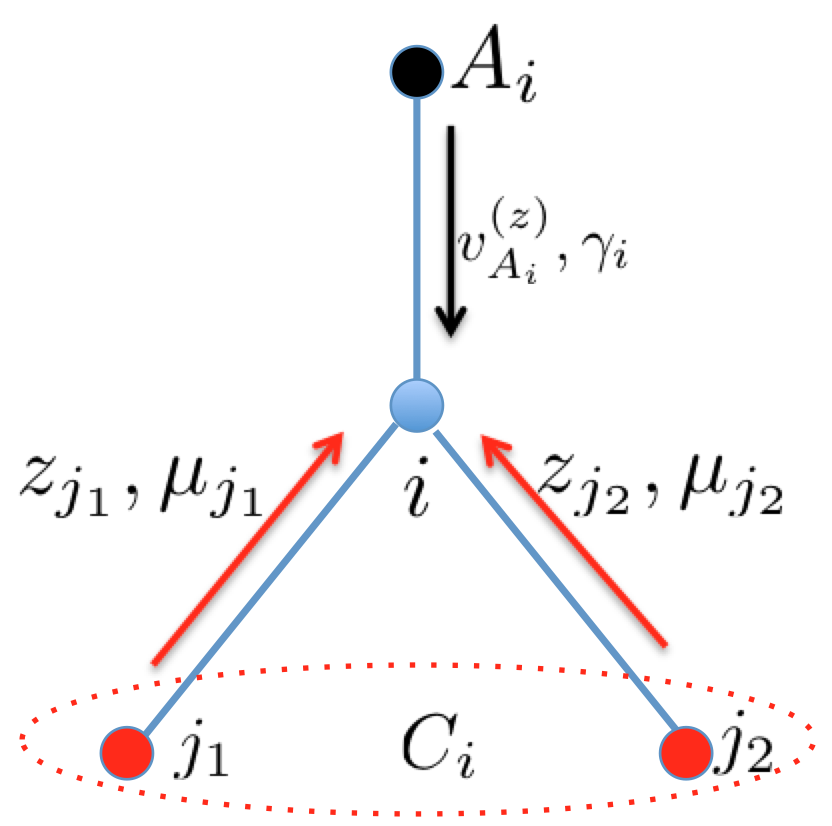}
	\label{fig:msg_x}
	}
\subfloat[z-update] {
	\includegraphics[scale=0.25]{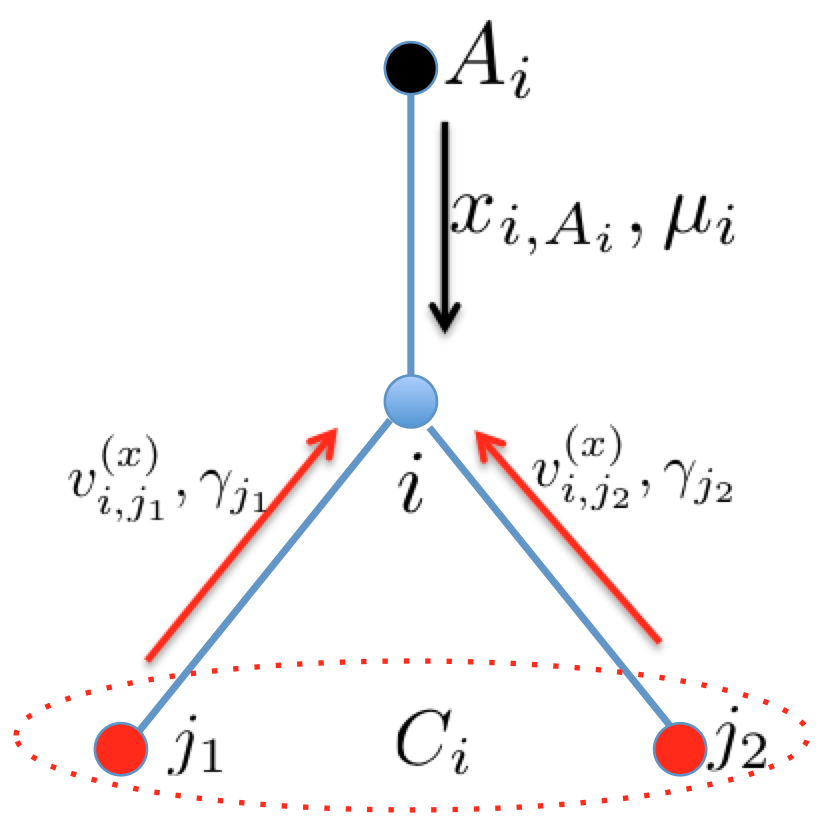}
	\label{fig:msg_z}
	}
\caption{Message exchange in the x and z-update step for agent $i$.}
\label{fig:msg_update}
\end{figure}

\begin{figure*}
{\scriptsize
\begin{subequations}\label{eq:AL}
\begin{align}
&L_{\rho}(x,z,\lambda,\gamma,\mu)\label{eq:AL1}\\
=&\sum_{i\in\mathN}\left(f_i(p_i^{(z)})+\langle\lambda_i,x_i-z_i\rangle+\langle\mu_i,x_{i,A_i}-z_i\rangle+\sum_{j\in C_i}\langle\gamma_j,v_{i,j}^{(x)}-v_{i}^{(z)}\rangle+\frac{\rho}{2}\left(\|x_i-z_i\|^2_s + \|x_{i,A_i}-z_i\|_n^2+\sum_{j\in C_i}\|v_{i,j}^{(x)}-v_{i}^{(z)}\|_2^2 \right)\right)\label{eq:AL2}\\
=&\sum_{i\in\mathN}\left(f_i(p_i^{(z)})+\langle\lambda_i,x_i-z_i\rangle+\sum_{j\in C_i}\langle\mu_j,x_{j,i}-z_j\rangle+\langle\gamma_i,v_{A_i,i}^{(x)}-v_{A_i}^{(z)}\rangle+\frac{\rho}{2}\left(\|x_i-z_i\|_s^2 + \sum_{j\in C_i}\|x_{j,i}-z_j\|_n^2+\|v_{A_i,i}^{(x)}-v_{A_i}^{(z)}\|_2^2 \right)\right)\label{eq:AL3}
\end{align}
\end{subequations}
}
\end{figure*}

Next, we demonstrate that the problem in \eqref{eq:eropf} can be solved in a distributed manner using ADMM, i.e. both the $x$-update \eqref{eq:xupdate} and $z$-update \eqref{eq:zupdate} can be decomposed into small subproblems that can be solved simultaneously by each agent $i$. For ease of presentation, we remove the iteration number $k$ in \eqref{eq:update} for all the variables, which will be updated accordingly after each subproblem is solved. The augmented Lagrangian for modified ROPF problem is given in \eqref{eq:AL}.
Note that in \eqref{eq:AL}, $x_{i,A_i}$ and $z_i$ consist of different components and $x_{i,A_i}-z_i$ is composed of the components that appear in both $x_{i,A_i}$ and $z_i$, \emph{i.e.} $x_{i,A_i}-z_i:=\left(S_{i,A_i}^{(x)}-S_{i}^{(z)},\ell_{i,A_i}^{(x)}-\ell_{i}^{(z)}\right)$.

In the $x$-update, each agent $i$ jointly solves the following $x$-update \eqref{eq:xupdate}.
\bq
\arg\min_{x\in\mathK_x}L_\rho(x,z,\lambda,\gamma,\mu)=\arg\min_{x\in\mathK_x}\sum_{i\in\mathN}G_i(x_i,v_{A_i,i}^{(x)},\{x_{j,i},j\in C_i\}),  \label{eq:xupdate_jointly}
\eq
where $G_i(x_i,v_{A_i,i}^{(x)},\{x_{j,i},j\in C_i\})$ is obtained from \eqref{eq:AL2} and
\begin{align*}
&G_i(x_i,v_{A_i,i}^{(x)},\{x_{j,i},j\in C_i\})\\
:=&\langle\lambda_i,x_i\rangle+\sum_{j\in C_i}\langle\mu_j,x_{j,i}\rangle+\langle\gamma_i,v_{A_i,i}^{(x)}\rangle+\frac{\rho}{2}\left(\|x_i-z_i\|_2^2 + \sum_{j\in C_i}\|x_{j,i}-z_j\|_2^2+\|v_{A_i,i}^{(x)}-v_{A_i}^{(z)}\|_2^2 \right)
\end{align*}
The corresponding subproblem for each agent $i$ that jointly solves \eqref{eq:xupdate_jointly} is
{
\begin{align}
\min \ &G_i(x_i,v_{A_i,i}^{(x)},\{x_{j,i},j\in C_i\})\nonumber\\
\mathrm{over}\  & x_i,v_{A_i,i}^{(x)},\{x_{j,i},j\in C_i\} \label{eq:xagent}\\
\mathrm{s.t.} \ & v^{(x)}_{A_i,i}-v^{(x)}_{i}+z_i\left(S^{(x)}_{i}\right)^*+S^{(x)}_{i}z_i^*-\ell^{(x)}_{i}|z_i|^2=0 \nonumber  \\
& \sum_{i\in C_i}\left(S^{(x)}_{j,i}-z_j\ell^{(x)}_{j,i}\right)-S^{(x)}_{i} +s^{(x)}_{i}=0 \nonumber
\end{align}
}

Prior to solving \eqref{eq:xagent}, each agent $i$ needs to collect $\left(v_{A_i}^{(z)},\gamma_i\right)$ from its parent and $(z_j,\mu_j)$ from all of its children $j\in C_i$.  The message exchanges in the x-update is illustrated in Fig. \ref{fig:msg_x}.

Next, we show how to solve \eqref{eq:xagent} in closed form. For each $i$, we can stack the real and imaginary part of the variables $(x_i,v_{A_i,i}^{(x)},\{x_{j,i},j\in C_i\})$ in a vector with appropriate dimensions and denote it as $\tilde x$. Then the subproblem solved by agent $i$ in the $x$-update \eqref{eq:xagent} takes the following form: 
\bq
\min_{\tilde x} \ \ \frac{1}{2}{\tilde x}^HA {\tilde x}+c^H{\tilde x} \quad \text{s.t. } B{\tilde x}=0  \label{eq:xupdate1}
\eq
where $A$ is a positive diagonal matrix, $B$ is a full row rank real matrix and $c$ is a real vector. $A,c,B$ are derived from \eqref{eq:xagent}. There exists a closed form expression for \eqref{eq:xupdate1} given by
\bqn
\tilde x=\left(A^{-1}B^H(BA^{-1}B^H)^{-1}BA^{-1}-A^{-1}\right)c
\eqn

In the $z$-update, each agent $i$ updates $z_i$ by jointly solving the z-update \eqref{eq:zupdate}
\bq
\arg\min_{z\in\mathK_z}L_\rho(x,z,\lambda,\gamma,\mu)=\arg\min_{z\in\mathK_z}\sum_{i\in\mathN}H_i(z_i),\label{eq:zupdate_jointly}
\eq
where $H_i(z_i)$ is obtained from \eqref{eq:AL3} and
\begin{align*}
H_i(z_i):=f_i(p_i^{(z)})-\langle\lambda_i,z_i\rangle-\langle\mu_i,z_i\rangle-\sum_{j\in C_i}\langle\gamma_j,v_{i}^{(z)}\rangle+\frac{\rho}{2}\left(\|x_i-z_i\|^2_2 + \|x_{i,A_i}-z_i\|_2^2+\sum_{j\in C_i}\|v_{i,j}^{(x)}-v_{i}^{(z)}\|_2^2\right)
\end{align*}
The corresponding subproblem for each agent $i$ that jointly solves \eqref{eq:zupdate_jointly} is 

\bq
\min && H_i(z_i)   \nonumber\\
\mathrm{over} && z_i \label{eq:zagent}\\
\mathrm{s.t.}&&    |S_i^{(z)}|^2\leq v_i^{(z)}\ell_i^{(z)} \nonumber   \\
&& s_i^{(z)}\in \mathI_i  \nonumber \\
&& \underline v_i\leq v_i^{(z)}\leq \overline v_i \nonumber
\eq

Prior to solving \eqref{eq:zagent}, each agent $i$ needs to collect $(x_{i,A_i},\mu_i)$ from its parent and $(v_{i,j}^{(x)},\gamma_j)$ from all of its children $j\in C_i$.  The message exchanges in the z-update is illustrated in Fig. \ref{fig:msg_z}.  

Next, we show how to solve \eqref{eq:zagent} in closed form. Note that 
{\small
\begin{align}
H_i(z_i):=&f_i(p_i^{(z)})-\langle\lambda_i,z_i\rangle-\langle\mu_i,z_i\rangle-\sum_{j\in C_i}\langle\gamma_j,v_{i}^{(z)}\rangle+\frac{\rho}{2}\left(\|x_i-z_i\|^2_2 + \|x_{i,A_i}-z_i\|_2^2+\sum_{j\in C_i}\|v_{i,j}^{(x)}-v_{i}^{(z)}\|_2^2\right)\nonumber\\
=&\rho\left(|S_i^{(z)}-\hat S_i|^2+|\ell_i^{(z)}-\hat \ell_i|^2+\frac{|C_i|+1}{2}|v_i^{(z)}-\hat v_i|^2\right)+f_i(p_i^{(z)})+\frac{\rho}{2}\|s_i^{(z)}-\hat s_i\|_2^2+\text{constant} \label{eq:zupdate_square}
\end{align} 
}
We use square completion to obtain \eqref{eq:zupdate_square} and the variables labeled with hat are some constants. Then \eqref{eq:zagent} can be furthered decomposed into two subproblems as below. The first one is 
\bq
\min && \left|S_i^{(z)}-\hat S_i\right|^2+\left|\ell_i^{(z)}-\hat \ell_i\right|^2+\frac{|C_i|+1}{2}\left|v_i^{(z)}-\hat v_i\right|^2\nonumber\\
\mathrm{over} && v_i^{(z)},\ell_i^{(z)}, S_i^{(z)} \label{eq:zagent1}\\
\mathrm{s.t.}&&     \left|S_i^{(z)}\right|^2\leq v_i^{(z)}\ell_i^{(z)} \nonumber\\
 && \underline v_i\leq v_i^{(z)}\leq \overline v_i \nonumber
\eq

The optimization problem in \eqref{eq:zagent1} has a quadratic objective, a second order cone constraint and a box constraint. We illustrate in Appendix \ref{app:newsolver} the procedure that solves \eqref{eq:zagent1}. Compared with using generic iterative solver, the procedure is computationally efficient since it only requires to solve the zero of three polynomials with degree less than or equal to $4$, which have closed form expression.

The second problem is 
\bq
\min && f_i\left(p_i^{(z)}\right)+\frac{\rho}{2}\left\|s_i^{(z)}-\hat s_i\right\|_2^2 \nonumber\\
\mathrm{over} &&  s_i^{(z)}  \label{eq:zagent2} \\
\mathrm{s.t.}&&  s_i^{(z)}\in \mathI_i \nonumber
\eq

Recall that $f_i\left(p_i^{(z)}\right):=\frac{\alpha_i}{2} p_i^2+\beta_i p_i$ as in \eqref{eq:objective}. If $\mathI_i$ takes the form of \eqref{eq:S2}, the closed form solution to \eqref{eq:zagent2} is 
\bqn
p_i^{(z)}=\left[\frac{\rho\hat p_i-\beta}{\alpha+\rho}\right]_{\underline p_i}^{\overline p_i} \quad
q_i^{(z)}=\left[\hat q_i\right]_{\underline q_i}^{\overline q_i}
\eqn
where $[x]_a^b:=\min\{a,\max\{x,b\}\}$. If $\mathI_i$ takes the form of \eqref{eq:S1}, there also exists a closed form expression to \eqref{eq:zagent2} and the procedure is relegate to Appendix \ref{app:solver2}.

Finally, we specify the initialization and stopping criteria for the algorithm. A good initialization usually reduce the number of iterations for convergence. We use the following initialization suggested by our empirical results. We first initialize the $z$ variable. The voltage magnitude square $v_i^{(z)}=1$. The power injection $s_i^{(z)}$ is picked up from a feasible point in the feasible region $\mathI_i$. The branch power $S_i^{(z)}$ is the aggregate power injection $s_i^{(z)}$ from the nodes connected by line $i$ (Note that the network has a tree topology.). The branch current $\ell_{i}^{(z)}=\frac{|S_i^{(z)}|^2}{v_i^{(z)}}$ according to \eqref{eq:bfm3}. The $x$ variables are initialized using the corresponding $z$ variable according to \eqref{eq:ndup}-\eqref{eq:sdup}. Intuitively, the above initialization procedure can be interpreted as a solution to the branch flow equation \eqref{eq:bfm} assuming zero impedance on all the lines. 

For the stopping criteria, there is no general rule for ADMM based algorithm and it usually hinges on the problem \cite{boyd2011distributed}. In \cite{boyd2011distributed}, it is suggested that a reasonable stopping criteria is that both the primal residual $r^k$ defined in \eqref{eq:pfeasible} and the dual residual $s^k$ defined in \eqref{eq:dfeasible} are within $10^{-3}\sqrt{|\mathN|}$ or $10^{-4}\sqrt{|\mathN|}$. The stopping criteria we adopt is that both $r^k$ and $s^k$ are below $10^{-4}\sqrt{|\mathN|}$ and empirical results show that the solution is accurate enough. The pseudo code for the algorithm is summarized in Table \ref{tab:alg}.

\begin{table}
\caption{Distributed algorithm of OPF}
\centering
\begin{tabular}{|p{8cm}|}
\hline
{\bf Distributed Algorithm of OPF}\\
\hline
{\bf Input}: network $\mathT$, power injection region $\mathI_i$, voltage region $(\underline v_i,\overline v_i)$,\\
\hspace{0.6cm} line impedance $z_i$ for $i\in\mathN$.\\
{\bf Output}: voltage $v$, power injection $s$.\\
\hline
1. Initialize the $x$ and $z$ variables. \\
2. Iterate the following step until both the primal residual $s^k$ \eqref{eq:pfeasible} and the dual residual $r^k$ \eqref{eq:dfeasible} are below $10^{-4}\sqrt{|\mathN|}$.\\
\hspace{0.2cm} a. In the $x$-update, each agent $i$ solves \eqref{eq:xagent} to update $x$.\\
\hspace{0.2cm} b. In the $z$-update, each agent $i$ solves \eqref{eq:zagent} to update $z$.\\
\hspace{0.2cm} c. In the multiplier update, update $\lambda,\mu,\gamma$ by  \eqref{eq:mupdate}.\\
\hline
\end{tabular}
\label{tab:alg}
\end{table}


\section{Case Study}\label{sec:case}

In this section, we first demonstrate the scalability of the distributed algorithm proposed in section \ref{sec:alg} by testing it on the model of a 2,065-bus distribution circuit in the service territory of the Southern California Edison (SCE). In particular, we also show the advantage of deriving closed form expression by comparing the computation time of solving the subproblems between off-the-shelf solver (CVX \cite{grant2008cvx}) and our algorithm. Second, we simulate the proposed algorithm on networks of different sizes to understand the factors that affect the convergence rate. The algorithm is implemented in Matlab 2014 and run on a Macbook pro 2014 with i5 dual core processor.

\subsection{Simulation on a 2,065 bus circuit}\label{sec:case1}

In the 2,065 bus distribution network, there are 1,409 household loads whose power consumptions are within 0.07kw--7.6kw and 142 commercial loads, whose power consumptions are within 5kw--36.5kw. There are 135 rooftop PV panels, whose nameplates are within 0.7--4.5kw, distributed across the 1,409 houses. 

The network is unbalanced three phase.
We assume that the three phases are decoupled such that the network becomes identical single phase network. The voltage magnitude at each load bus is allowed to lie within 
$[0.95,1.05]$ per unit (pu), i.e. $\overline v_i=1.05^2$ and $\underline v_i=0.95^2$ for $i\in\mathN_+$. The control devices are the rooftop PV panels whose real and reactive power injections are controlled. The objective is to minimize power loss across the network, namely $\alpha_i=0,\beta_i=1$ for $i\in\mathN$, where $\alpha_i,\beta_i$ are coefficients in the objective function and defined in \eqref{eq:objective}. Each bus is an agent and there are 2,065 agents in the network 
that solve the OPF problem in a distributed manner.

\begin{figure}
\centering
\subfloat[Primal and dual residual] {
	\includegraphics[scale=0.3]{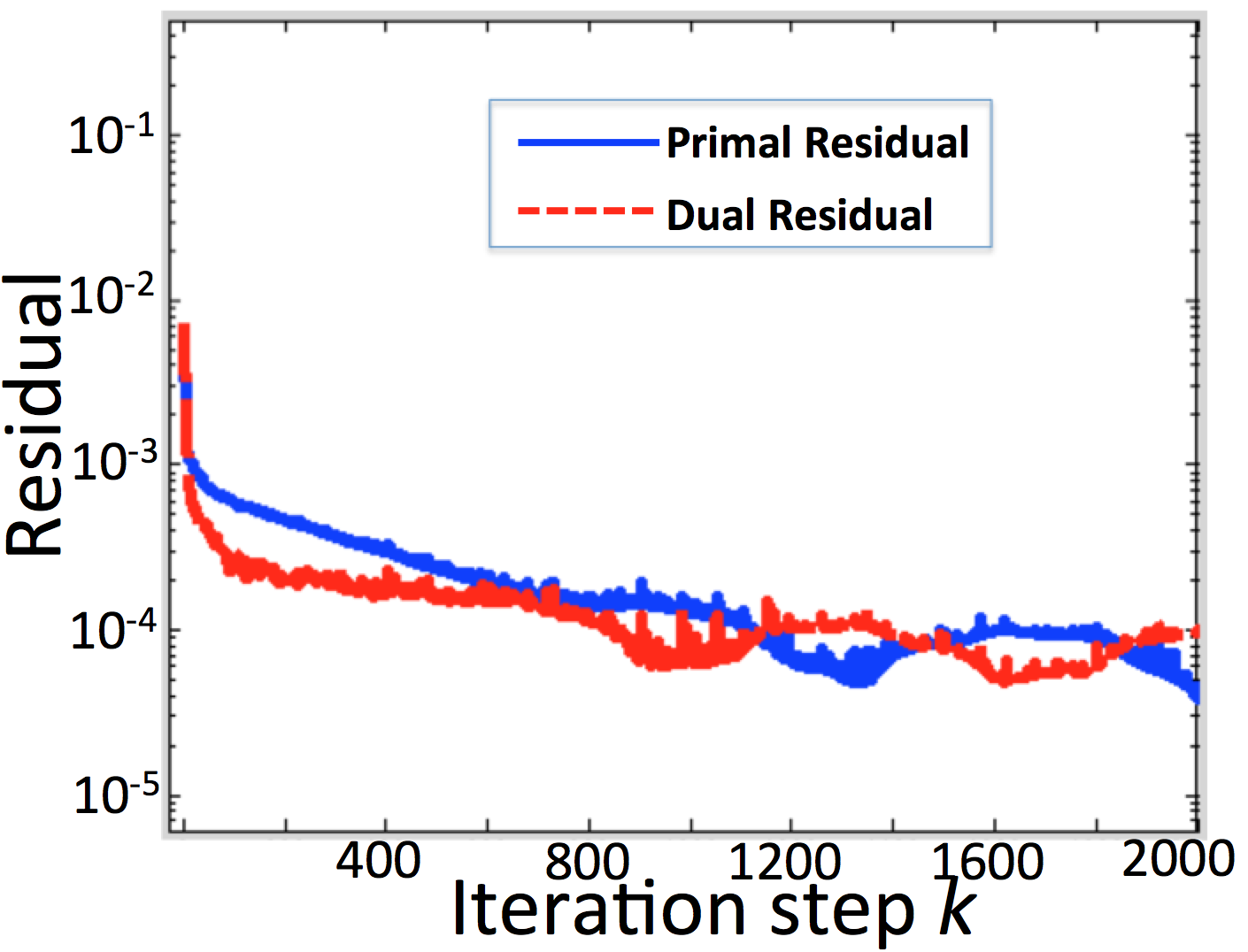}
	\label{fig:residual}
	}
\subfloat[Objective value]{
	\includegraphics[scale=0.3]{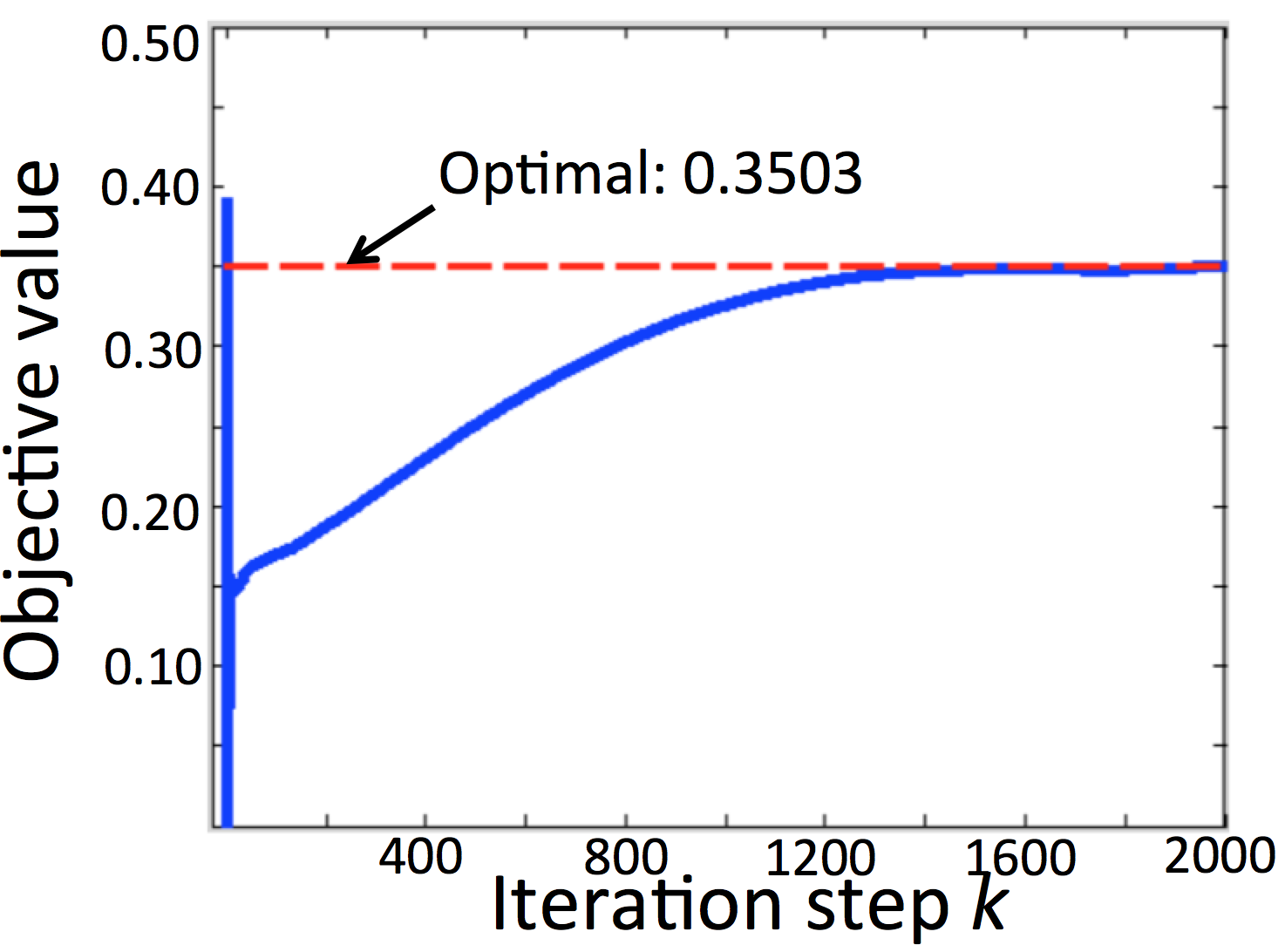}
	\label{fig:objective}
	}
\caption{Simulation results for 2065 bus distribution network.}
\label{fig:rossi2065}
\end{figure}

We mainly focus on the time of convergence (ToC) for
the proposed distributed algorithm. The algorithm is run on a single machine. To roughly estimate the ToC (excluding communication overhead) if the algorithm is run on distributed machines, we divide the total time by the number of agents. Recall that the stopping criteria is that both the primal and dual residual are below $10^{-4}\sqrt{|\mathN|}$ and Figure \ref{fig:residual} illustrates the evolution of $r^k/\sqrt{|\mathN|}$ and $s^k/\sqrt{|\mathN|}$ versus iterations $k$. The stopping criteria are satisfied after $1,114$ iterations. The evolution of the objective value is illustrated in Figure \ref{fig:objective}. It takes 1,153s to run 1,114 iterations on a single computer.
Then the ToC is roughly 0.56s if we implement the algorithm in a distributed manner not counting communication overhead. 

Finally, we show the advantage of closed form solution by comparing the computation time of solving the subproblems by an off-the-shelf solver (CVX) and by our algorithm. In particular, we compare the average computation time of solving the subproblem in both the $x$-update and the $z$-update step. In the $x$-update step, the average time required to solve the subproblem is $1.7\times10^{-4}$s for the proposed algorithm but $0.2$s for CVX.  In the $z$-update step, the average time required to solve the subproblem is $5.1\times10^{-4}$s for the proposed algorithm but $0.3$s for CVX. Thus, each ADMM iteration only takes about $6.8\times 10^{-4}$s for the proposed algorithm but $0.5$s for using iterative algorithm, which is a 1,000x speedup. 

\subsection{Rate of Convergence}

\begin{table}
\caption{Statistics of different networks}
\begin{center}
\begin{tabular}{|c|c|c|c|c|}
\hline
Network & Diameter & Iteration & Total Time(s) & Avg time(s)\\
\hline
2065Bus& 64 & 1114 &1153 & 0.56\\
\hline
1313Bus& 53 &671 &471 & 0.36\\ 
\hline
792Bus& 45 &524 &226 & 0.29\\ 
\hline
363Bus& 36 &289 &112 & 0.24\\ 
\hline
108Bus& 16 & 267 &16 & 0.14\\ 
\hline
\end{tabular}
\end{center}
\label{tab:statistics}
\end{table}

In section \ref{sec:case1}, we demonstrate that the proposed distributed algorithm can dramatically reduce the computation time within each iteration and therefore is scalable to a large practical 2,065 bus distribution network. The time of convergence(ToC) is determined by both the computation time required within each iteration and the number of iterations. In this subsection, we study the number of iterations, namely rate of convergence.

To our best knowledge, most of the works on convergence rate for ADMM based algorithms study how the primal/dual residual changes as the number of iterations increases. Specifically, it is proved  in \cite{wei20131,hong2012linear} that the general ADMM based algorithms converge linearly under certain assumptions. Here, we consider the rate of convergence from another two factors, network size $N$ and diameter $D$, \emph{i.e.} given the termination criteria in Table \ref{tab:alg}, the impact from network size and diameter on the number of iterations. The impacts from other factors, \emph{e.g.} form of objective function and constraints, etc. are beyond the scope of this paper. 

First, we simulate the algorithm on different networks (that are subnetworks of the 2,065-bus system) and some statistics are given in Table \ref{tab:statistics}. For simplicity, we assume the number of iterations $T$ to converge takes the linear form $T=a N+b D$. Using the data in Table \ref{tab:statistics}, the parameters $a=0.34$, $b=5.53$ give the least square error. It means that the network diameter has a stronger impact than the network size on the rate of convergence.

\begin{figure}
\centering
\subfloat[Line network] {
	\includegraphics[scale=0.6]{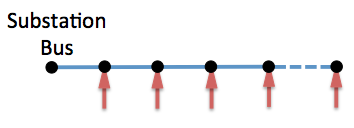}
	\label{fig:line}
	}
\subfloat[Fat tree network]{
	\includegraphics[scale=0.6]{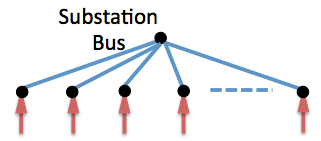}
	\label{fig:fattree}
	}
\caption{Topology for tree and fat tree networks.}
\label{fig:network}
\end{figure}

 To further illustrate the phenomenon, we simulate the algorithm on two extreme cases: 1) Line network in Fig. \ref{fig:line}, whose diameter is the largest given the network size, 2) Fat tree network in Fig. \ref{fig:fattree}, whose diameter is the smallest ($2$) given the network size. In Table \ref{tab:statistics1}, we record the number of iterations for both line and fat tree network of different sizes. For line network, the number of iterations increases notably as the size increases. For fat tree network, the trend is less obvious compared to line network.

\begin{table}
\caption{Statistics of line and fat tree networks}
\begin{center}
\begin{tabular}{|c|c|c|}
\hline
Size & $\#$ of iterations (Line) & $\#$ of iterations (Fat tree) \\
\hline
$5$ & $43$ & $31$ \\
\hline
$10$ & $123$  & $51$\\
\hline
$15$ & $198$ & $148$\\
\hline
$20$ & $286$ & $87$\\
\hline 
$25$ & $408$ & $173$\\
\hline 
$30$ & $838$ & $119$\\
\hline
$35$ & $1471$ & $187$\\ 
\hline
$40$ & $2201$ & $109$\\
\hline
$45$ & $2586$ & $182$\\
\hline
$50$ & $3070$ & $234$\\
\hline
\end{tabular}
\end{center}
\label{tab:statistics1}
\end{table}

\section{Conclusion}\label{sec:conclusion}
In this paper, we have developed a distributed algorithm for optimal power flow problem based on alternating direction method of multiplier for balanced radial distribution network. We have derived a closed form solution for the subproblems solved by each agent thus significantly reducing the computation time. Preliminary simulation shows that the algorithm is scalable to a 2,065-bus system and the optimization subproblem in each ADMM iteration is solved 1,000x faster than generic optimization solver.

\bibliographystyle{IEEEtran}
\bibliography{texcode/PowerRef-201202,texcode/reference}{}

\appendices
\section{}\label{app:newsolver}

Denote $z_1:=\mathbf{Re}(S_i^{(z)})$, $z_2:=\mathbf{Im}(S_i^{(z)})$, $z_3:=\sqrt{\frac{|C_i|+1}{2}}v_i^{(z)}$ and $z_4:=\ell_i^{(z)}$. Then the optimization problem \eqref{eq:zagent1} can be written equivalently as 
\begin{subequations} \label{eq:zagent1_2}
\bq
\min && \sum_{i=1}^4(z_i^2+c_iz_i) \label{eq:obj1}\\
\mathrm{over} && z  \nonumber \\
\mathrm{s.t. } && \frac{z_1^2+z_2^2}{z_3}\leq k^2 z_4   \label{eq:con11}\\
&& \underline z_3 \leq z_3\leq \overline z_3  \label{eq:con12}
\eq
\end{subequations}
where $\overline z_3> \underline z_3> 0$ and $c_i$, $k$ are constants that hinges on the constants in \eqref{eq:zagent1}.

Below we will derive a procedure that solves \eqref{eq:zagent1_2}. Let $\mu\geq 0$ denote the Lagrangian multiplier for constraint \eqref{eq:con11} and $\underline \lambda,\overline \lambda\geq 0$ denote the Lagrangian multipliers for constraint \eqref{eq:con12}, then the Lagrangian of P1 is 
\bqn
L(z,\mu,\lambda)=\sum_{i=1}^4(z_i^2+c_iz_i)+\mu\left( \frac{z_1^2+z_2^2}{z_3}-k^2z_4\right)+ \overline \lambda(z_3-\overline z_3)-\underline \lambda(z_3-\underline z_3)
\eqn

The KKT optimality conditions imply that the optimal solution $z^*$ together with the multipliers $\mu^*,\underline \lambda^*,\overline \lambda^*$ satisfy the following equations. For ease of notations, we sometimes skip the superscript $\star$ of the variables in the following analysis.
\begin{subequations}\label{eq:kkt1}
\begin{align}
&2z_1+c_1+2\mu\frac{z_1}{z_3}=0\\
&2z_2+c_2+2\mu\frac{z_2}{z_3}=0\\
&2z_3+c_3-\mu\frac{z_1^2+z_2^2}{z_3^2}+\overline \lambda-\underline \lambda=0\\
&2z_4+c_4-k^2\mu=0\\
& \overline \lambda(z_3-\overline z_3)=0 \quad \overline \lambda\geq 0 \quad z_3\leq \overline z_3 \\
& \underline \lambda(\underline z_3-z_3)=0 \quad \underline \lambda\geq 0 \quad z_3\geq \underline z_3 \\
& \mu\left( \frac{z_1^2+z_2^2}{z_3}-k^2z_4\right)=0 \quad \mu\geq 0 \quad \frac{z_1^2+z_2^2}{z_3}\leq k^2z_4 
\end{align}
\end{subequations}

\begin{lemma}\label{lem:p1unique}
There exists a unique solution $(z^*,\mu^*,\underline \lambda^*,\overline \lambda^*)$ to \eqref{eq:kkt1} if $\overline z_3> \underline z_3\geq 0$.
\end{lemma}
\begin{IEEEproof}
P1 is feasible since $z=(0,0,\underline z_3,1)$ satisfies  \eqref{eq:con11}-\eqref{eq:con12}. In addition, P1 is a strictly convex optimization problem since the objective \eqref{eq:obj1} is a strictly convex function of $z$ and the constraints \eqref{eq:con11} and \eqref{eq:con12} are also convex. Hence, there exists a unique solution $z^*$ to P1, which indicates there exists a unique solution $(z^*,\mu^*,\underline \lambda^*,\overline \lambda^*)$ to the KKT optimality conditions \eqref{eq:kkt1}.
\end{IEEEproof}

Lemma \ref{lem:p1unique} says that there exists a unique solution to \eqref{eq:kkt1}, which is also the optima to P1. In the following, we will solve \eqref{eq:kkt1} through enumerating value of the multipliers $\mu,\overline\lambda,\underline\lambda$. Specifically, we first assume $\mu^*=0$ (Case 1 below), which is equivalent to assume constraint \eqref{eq:con11} is inactive. If there is a feasible solution to \eqref{eq:kkt1}, it is the unique solution to \eqref{eq:kkt1}. Otherwise, we assume $\mu^*=0$ (Case 2 below), which is equivalent to assume that the equality is obtained at optimality in \eqref{eq:con11}.

\flushleft\textbf{Case 1:} If $\mu=0$, \eqref{eq:kkt1} becomes
\begin{subequations}\label{eq:case1}
\begin{align}
&2z_1+c_1=0 \label{eq:case11}\\
&2z_2+c_2=0\\
&2z_3+c_3+\overline \lambda-\underline \lambda=0\\
&2z_4+c_4=0 \label{eq:case14}\\
& \overline \lambda(z_3-\overline z_3)=0 \quad \overline \lambda\geq 0 \quad z_3\leq \overline z_3 \label{eq:case15}\\
& \underline \lambda(z_3-\underline z_3)=0 \quad \underline \lambda\geq 0 \quad z_3\geq \underline z_3  \label{eq:case16}\\
& \frac{z_1^2+z_2^2}{z_3}\leq k^2z_4 \label{eq:case17}
\end{align}
\end{subequations}
The solution to \eqref{eq:case11}--\eqref{eq:case16} ignoring \eqref{eq:case17} is 
\begin{align*}
&z_1=-\frac{c_1}{2}, \ z_2=-\frac{c_2}{2}, \ z_3=\left[-\frac{c_3}{2}\right]_{\underline z_3}^{\overline z_3},  \ z_4=-\frac{c_4}{2}\\
& \overline \lambda=-(2z_3+c_3)1_{\{z_3=\overline z_3\}}, \underline \lambda=-(2z_3+c_3)1_{\{z_3=\underline z_3\}} 
\end{align*}
and if the solution satisfies \eqref{eq:case17}, it is the solution to \eqref{eq:kkt1}. Otherwise, we go to Case 2.

\flushleft\textbf{Case 2:} If $\mu>0$, \eqref{eq:kkt1} becomes
\begin{subequations}\label{eq:kkt2}
\begin{align}
&2z_1+c_1+2\mu\frac{z_1}{z_3}=0\label{eq:kkt21}\\
&2z_2+c_2+2\mu\frac{z_2}{z_3}=0\label{eq:kkt22}\\
&2z_3+c_3-\mu\frac{z_1^2+z_2^2}{z_3^2}+\overline \lambda-\underline \lambda=0\label{eq:kkt23}\\
& \overline \lambda(z_3-\overline z_3)=0 \quad \overline \lambda\geq 0 \quad z_3\leq \overline z_3 \\
& \underline \lambda(z_3-\underline z_3)=0 \quad \underline \lambda\geq 0 \quad z_3\geq \underline z_3 \label{eq:kkt25}\\
&\mu=\frac{1}{k^2}(2z_4+c_4)\label{eq:kkt26}\\
&z_4=\frac{z_1^2+z_2^2}{k^2z_3} \label{eq:kkt27}
\end{align}
\end{subequations}
Substitute \eqref{eq:kkt27} into \eqref{eq:kkt26}, we obtain
\bq\label{eq:kkt28}
\mu=\frac{1}{k^2}(2z_4+c_4)=\frac{2(z_1^2+z_2^2)}{k^4z_3}+\frac{c_4}{k^2}
\eq
Then substituting \eqref{eq:kkt26} into \eqref{eq:kkt21}-\eqref{eq:kkt25}, we can write \eqref{eq:kkt2}  equivalently as

\begin{subequations}\label{eq:kkt3}
\begin{align}
&2+\frac{c_1}{z_1}+4\frac{(z_1^2+z_2^2)}{k^4z_3^2}+\frac{2c_4}{k^2z_3}=0 \label{eq:kkt31}\\
&2+\frac{c_2}{z_2}+4\frac{(z_1^2+z_2^2)}{k^4z_3^2}+\frac{2c_4}{k^2z_3}=0\label{eq:kkt32}\\
&2+\frac{c_3}{z_3}-2\frac{(z_1^2+z_2^2)^2}{k^4z_3^4}-c_4\frac{z_1^2+z_2^2}{z_3^3}+\frac{\overline \lambda-\underline \lambda}{z_3}=0\label{eq:kkt33}\\
& \overline \lambda(z_3-\overline z_3)=0 \quad \overline \lambda\geq 0 \quad z_3\leq \overline z_3 \\
& \underline \lambda(z_3-\underline z_3)=0 \quad \underline \lambda\geq 0 \quad z_3\geq \underline z_3 
\end{align}
\end{subequations}
where \eqref{eq:kkt31}--\eqref{eq:kkt33} are obtained through dividing both sides of \eqref{eq:kkt21}--\eqref{eq:kkt23} by $z_1$, $z_2$ and $z_3$, respectively. The variables $\mu,z_4$ can be recovered via \eqref{eq:kkt26} and \eqref{eq:kkt27} after we solve \eqref{eq:kkt3}.

By \eqref{eq:kkt31} and \eqref{eq:kkt32}, 
\bqn
\frac{c_1}{z_1}=\frac{c_2}{z_2}
\eqn
Denote $p:=\frac{z_1}{c_1z_3}=\frac{z_2}{c_2z_3}$. Then \eqref{eq:kkt3} is equivalent to the following equations. 
\begin{subequations}\label{eq:kkt4}
\begin{align}
&p=\frac{z_1}{c_1z_3}=\frac{z_2}{c_2z_3} \label{eq:kkt41}\\
&2+\frac{1}{pz_3}=-\left(\frac{4(c_1^2+c_2^2)}{k^4}p^2+2\frac{c_4}{k^2z_3}\right) \label{eq:kkt42}\\
&2+\frac{c_3}{z_3}=\frac{2(c_1^2+c_2^2)^2}{k^4}p^4+\frac{c_4(c_1^2+c_2^2)}{k^2}\frac{p^2}{z_3}+\frac{\underline \lambda-\overline \lambda}{z_3}\label{eq:kkt43}\\
& \overline \lambda(z_3-\overline z_3)=0 \quad \overline \lambda\geq 0 \quad z_3\leq \overline z_3 \\
& \underline \lambda(z_3-\underline z_3)=0 \quad \underline \lambda\geq 0 \quad z_3\geq \underline z_3 
\end{align}
\end{subequations}
where $\eqref{eq:kkt42}$ is obtained by substitute $z_2=c_2pz_3$ into \eqref{eq:kkt32}, \eqref{eq:kkt43} is obtained by substitute $z_1=c_1pz_3$ and $z_2=c_2pz_3$ into \eqref{eq:kkt33}. To solve \eqref{eq:kkt4}, we further divide our analysis into two sub-cases depending on whether $z_3^*$ hits the lower or upper bound.

\bi
\item {\bf Case 2.1}: $z_3^*=\overline z_3$. $(\underline \lambda=0,\overline\lambda>0)$ ( The case of $z_3^*=\underline z_3$ can be solved using similar procedure.) \\
We first substitute $z_3=\overline z_3$ into \eqref{eq:kkt42} and  have
\bq
\frac{4(c_1^2+c_2^2)}{k^4}p^3+\left(2\frac{c_4}{k^2\overline z_3}+2\right)p+\frac{1}{\overline z_3}=0,
\eq
whose solution\footnote{\label{footnote:1}There are potentially multiple solutions and we need to check all the real solution $p^*$ using the following procedure.} is denoted by $p^*$. Then substitute $p^*$ and $\overline z_3$ into \eqref{eq:kkt41}, we can recover $z_1^*$ and $z_2^*$. Then we can obtain $\mu^*,\overline \lambda^*$ using \eqref{eq:kkt28} and \eqref{eq:kkt33} by substituting $z_1^*,\ldots,z_4^*$. If $\mu^*,\overline \lambda^*\geq0$, they collectively solve \eqref{eq:kkt1}. Otherwise, we go to Case 2.2.

\item {\bf Case 2.2}: $\underline z_3<z_3^*<\overline z_3$ $(\underline \lambda,\overline\lambda=0)$.\\
Since $\overline \lambda$ and $\underline \lambda=0$, \eqref{eq:kkt4} reduces to 
\begin{subequations}\label{eq:kkt5}
\begin{align}
&p=\frac{z_1}{c_1z_3}=\frac{z_2}{c_2z_3} \label{eq:kkt51}\\
&2+\frac{1}{pz_3}=-\left(\frac{4(c_1^2+c_2^2)}{k^4}p^2+2\frac{c_4}{k^2z_3}\right) \label{eq:kkt52}\\
&2+\frac{c_3}{z_3}=\frac{2(c_1^2+c_2^2)^2}{k^4}p^4+\frac{c_4(c_1^2+c_2^2)}{k^2}\frac{p^2}{z_3} \label{eq:kkt53}
\end{align}
\end{subequations}
\ei
Dividing each side of \eqref{eq:kkt52} by \eqref{eq:kkt53} gives
\bqn
\frac{2z_3+\frac{1}{p}}{2z_3+c_3}=-\frac{2}{(c_1^2+c_2^2)p^2},
\eqn
which implies
\bq\label{eq:case2}
z_3=-\frac{(c_1^2+c_2^2)p+2c_3}{2((c_1^2+c_2^2)p^2+2)}
\eq
Then substitute $\eqref{eq:case2}$ into \eqref{eq:kkt52}, we have 
{\small
\bqn
\frac{(c_1^2+c_2^2)p^2+2}{(c_1^2+c_2^2)p^2+2c_3p}-\frac{2(c_1^2+c_2^2)}{k^4}p^2+\frac{2c_4((c_1^2+c_2^2)p^2+2)}{k^2((c_1^2+c_2^2)p+2c_3)}=1
\eqn
}
which is equivalent to 
{\small
\bqn
\frac{(c_1^2+c_2^2)^2}{k^4}p^4+\frac{c_1^2+c_2^2}{k^2}\left(\frac{2c_3}{k^2}-c_4\right)p^3+\left(c_3-\frac{2c_4}{k^2}\right)p-1=0
\eqn
}
whose solution is denoted by $p^*$. Substitute $p^*$ into \eqref{eq:case2}, we can recover $z_3^*$, then $z_1^*,z_2^*$ can be recovered via \eqref{eq:kkt51}. $\mu^*$ is recovered using \eqref{eq:kkt28}. If $\mu^*\geq 0$, the corresponding solution solves \eqref{eq:kkt1}.



\section{}\label{app:solver2}

If $\mathI_i$ takes the form of \eqref{eq:S1}, the optimization problem \eqref{eq:zagent2} takes the following form
\begin{subequations}\label{eq:app2}
\bq
\min_{p,q} && \frac{a_1}{2}p^2+b_1p+\frac{a_2}{2}q^2+b_2q\\
\text{s.t. } && p^2+q^2\leq c^2 \label{eq:app2:1}\\
&& p \geq 0 \label{eq:app2:2}
\eq
\end{subequations}
where $a_1,a_2,c>0$ $,b_1,b_2$ are constants. The solutions to \eqref{eq:app2} are given as below.
{\flushleft\bf Case 1}: $b_1\geq 0$.
\bqn
p^*=0 \qquad q^*=\left[-\frac{b_2}{a_2}\right]_{-c}^{c}
\eqn
{\flushleft\bf Case 2}: $b_1< 0$ and $\frac{b^2_1}{a^2_1}+\frac{b^2_2}{a_2^2}\leq c^2$.
\bqn
p^*=-\frac{b_1}{a_1} \qquad q^*=-\frac{b_2}{a_2}
\eqn
{\flushleft\bf Case 3}: $b_1< 0$ and $\frac{b^2_1}{a^2_1}+\frac{b^2_2}{a_2^2}> c^2$.\\
First solve the following equation in terms of variable $\lambda$:
\bq
b_1^2(a_2+2\lambda)^2+b_2^2(a_1+2\lambda)^2=(a_1+2\lambda)^2(a_2+2\lambda)^2 \label{eq:app3:1}
\eq
which is a polynomial with degree of $4$ and has closed form expression. There are four solutions to \eqref{eq:app3:1}, but there is only one strictly positive $\lambda^*$, which can be proved via the KKT conditions of \eqref{eq:app2}. Then we can recover $p^*,q^*$ from $\lambda^*$ using the following equations.
\bqn
 p^*=-\frac{b_1}{a_1+2\lambda^*} \quad \text{ and }\quad q^*=-\frac{b_2}{a_2+2\lambda^*} 
\eqn

The above procedure to solve \eqref{eq:app2} is derived from standard applications of the KKT conditions of \eqref{eq:app2}. For brevity, we skip the proof here.



\end{document}